\documentclass[nummat]{svjour}

\usepackage{amssymb}
\usepackage{amsmath}
\usepackage[final]{graphicx}
\usepackage{psfrag}
\usepackage{epsfig}
\usepackage{latexsym}
\usepackage{exscale}

\font\Bbb=msbm10 %at 12pt

\newcommand{\Rset}{\mbox{\Bbb R}}

\newcommand{\bR}{\mathbb{R}}
\newcommand{\bC}{\mathbb{C}}

\newcommand{\real}{{\rm Re}\,}
\newcommand{\imag}{{\rm Im}\,}
\newcommand{\weg}[1]{ }

\def \dt    { h }

\def \1{1\!{\rm l}}
\def \A{{\cal O}\!\!\iota}

\def \~{\widetilde}
\def\Log{\log}

\def\Lfunc{{\phi}}
\def\five{ }

\begin{document}

\title{Fast Runge-Kutta approximation of\\
inhomogeneous parabolic %differential 
equations}

\author{Mar\'\i a L\'opez-Fern\'andez$^1$, 
Christian Lubich$^2$,\\ 
Cesar Palencia$^1$, 
and Achim Sch\"adle$^3$}

\authorrunning{M.~L\'opez-Fern\'andez, C.~Lubich, C.~Palencia, and A.~Sch\"adle}

\institute{$^1$~Departamento de Matem\'atica Aplicada y Computaci\'on, 
Universidad de Valladolid, Valladolid, Spain.
~E-mail: {\tt \{marial, palencia\}@mac.cie.uva.es}\\
$^2$~Mathematisches Institut,
  Universit\"at T\"ubingen,
  Auf der Morgenstelle 10,\\
  D--72076 T\"ubingen,
  Germany. ~E-mail:  {\tt lubich@na.uni-tuebingen.de}\\
$^3$~ZIB Berlin, Takustr.~7, D-14195 Berlin, Germany.
~E-mail: {\tt schaedle@zib.de}     }

\date{\today}

\maketitle

\begin{abstract}  The result after $N$  steps of an
implicit Runge-Kutta time discretization of an inhomogeneous linear
parabolic differential equation is computed, up to accuracy $\varepsilon$, 
by solving only 
$$O\Big(\log N\, \log  \frac1\varepsilon \Big)
$$ 
linear systems
of equations. We derive, analyse, and numerically illustrate this
fast algorithm.
\end{abstract}


\section{Introduction}
In the method of lines, semi-discretization in space turns a
linear parabolic differential equation into a large, stiff system of
ordinary differential equations
\begin{equation}\label{ivp}
u'(t) + Au(t) = g(t), \quad\ u(0)=u_0,
\end{equation}
possibly with a mass matrix multiplying the time derivative.
This system is subsequently discretized in time, 
e.g., by the
implicit Euler method  with  step size $\dt$,
$$
(I+\dt A) u_n = u_{n-1} + \dt g(t_n),   \qquad n=1,\dots, N.
$$
The approximation $u_N$ for a prescribed step number $N$ is thus obtained by
solving a sequence of $N$ linear systems with a matrix of the form 
$\lambda +A$,
where $\lambda = 1/\dt$ in the implicit Euler method. 
For $N$~steps with a higher-order, $m$-stage Runge-Kutta method,
there are $mN$ such
linear systems, possibly with complex~$\lambda$ 
as in the 
excellent Radau~IIA methods.
Even if fast techniques such as multi-grid methods are used,
solving the linear systems of equations 
typically constitutes the main computational
cost, in particular for problems in complicated spatial geometries. 

In this paper we propose an algorithm to compute 
the implicit Runge-Kutta approximation $u_N$ at a fixed time $T=N\dt$, up to an
arbitrary accuracy $\varepsilon$, by doing
$N$ Runge-Kutta steps for
differential equations of the form
$y'(t) = \lambda y (t)+ g(t)$, each step
in parallel for $O(\Log (1/\varepsilon))$ 
complex parameters $\lambda$, and by
solving only 
$$
\hbox{$O(\log N \,\Log\frac1\varepsilon)$ linear systems}
$$
%(where $\Log x = \log x \cdot \log\log x$)
 with matrices
of the form $\lambda +A$, all of which can be solved in parallel.
The constant in this work estimate is moderate: for a
relative accuracy of $10^{-5}$ and $N\le 10^5$ time steps we need to solve
less than 100 linear systems!
For large step numbers $N$, the number of linear systems
is thus dramatically reduced, both in a sequential 
and in a parallel computational setting. 

The algorithm is highly efficient
for computing Runge-Kutta approximations to the solution of (\ref{ivp}) 
at a relatively small number of selected time points or of short subintervals, 
but it is not useful for computing {\it all} values $u_1,\dots, u_N$.

\pagebreak[2]
Basic ingredients of the algorithm are the following:
\begin{itemize}
\item the discrete variation-of-constants
formula for the Runge-Kutta method; 
\item the Cauchy integral representation of the approximations to
the operator exponential; 
\item the discretization of  the contour integrals, using
$O(\log N)$ contours with $O(\Log (1/\varepsilon))$ quadrature points each;
\item the discrete semigroup property, which permits us to reinterpret the
split sums as Runge-Kutta approximations to solutions of equations of the form
$y'(t)=\lambda y(t)+ g(t)$.
\end{itemize}

The algorithm given here is closely related to the fast convolution
algorithms developed in \cite{LuS02,SchLL}. 
The error analysis for the discretized contour integrals
follows the analysis of inverse Laplace transform approximations
in \cite{LoP04}.
 
Discretized contour integrals have been used previously
in several instances
in the numerical solution of parabolic equations: for homogeneous
problems ($g\equiv 0$) in \cite{ShST00} similarly to
Talbot's method \cite{Tal79} for the inversion of the Laplace transform
$(s+A)^{-1}u_0$, and more recently for
inhomogeneous problems \cite{ShST03,GaM} using the Laplace
transform of the inhomogeneity $g$ or assuming special properties,
in particular analyticity, of $g$. In contrast, the present
algorithm works directly with the discrete values $g(t)$ that are 
used in the Runge-Kutta discretization of (\ref{ivp}).
No smoothness conditions for $g$ are needed.
This is because the algorithm approximates 
the {\it discrete} result of the Runge-Kutta method, 
with an error that does not depend
on the smoothness of either the inhomogeneity or the solution. Of course,
to make sense,
the Runge-Kutta discretization of (\ref{ivp}) with the considered
step size $\dt$ should be sufficiently
accurate, which in turn does depend on the smoothness of $g$
(see \cite{LuO93} for Runge-Kutta error bounds for parabolic
equations in terms of
the data). 

About the differential 
equation (\ref{ivp}) we assume that $A$ is {\it sectorial}:
there exist real constants $M$ and $\sigma$ and an angle $\varphi<\frac\pi 2$ 
such that the
resolvent is bounded by
\begin{equation}\label{sector}
\| (\lambda + A)^{-1} \| \le \frac{M }{ |\lambda-\sigma|},
\quad\ \hbox{ for }\quad |\arg(\lambda -\sigma)|\le \pi-\varphi.
\end{equation}
Here $\|\cdot\|$ is the operator norm corresponding to a vector norm,
also denoted by $\|\cdot\|$. Clearly, for a symmetric positive semi-definite matrix $A$
the bound (\ref{sector}) holds in the Euclidean norm with 
$\sigma=0$ and $M=1/\sin\varphi$ for any positive angle $\varphi$.
More generally, condition (\ref{sector}) includes also non-symmetric
operators such as those arising in convection-diffusion equations.
In many situations resolvent bounds (\ref{sector}) 
in $L^p$~norms are
known to be inherited from the continuous problem by 
finite differences or finite elements, uniformly in the 
spatial discretization parameter 
(see, e.g., \cite{AsS94,BaTW02}). 

In Section~\ref{sec.dvc} we review the discrete variation-of-constants formula for
implicit Runge-Kutta methods, and in Section~\ref{Sec.Contour} 
we describe the discretization
of the contour integrals for the rational
approximations to the matrix exponential.
The fast algorithm is given in Section~\ref{Sec.FA}, including an extension to 
systems with a mass matrix.  A numerical example illustrates the
performance of the algorithm in Section~\ref{sec:NumEx}. Finally, 
Section~\ref{sec.error} analyses the error of the
contour integral discretization, which is the only error source in the
algorithm.

\section{The discrete variation-of-constants formula}
\label{sec.dvc}
In this preparatory section we recall
the discrete variation-of-constants
formula for implicit Runge-Kutta methods; cf., e.g.,~\cite{BrCT82}.

An implicit $m$-stage 
Runge-Kutta method applied to (\ref{ivp}) yields, at $t_n=n\dt$,
an approximation $u_n$ to $u(t_n)$, given recursively by
\begin{eqnarray}
v_{ni} &=& u_n + \dt\sum_{j=1}^m a_{ij} \Bigl( -A v_{nj} + g(t_n+c_j \dt) \Bigr),
\quad 1 \le i \le m,
\\
u_{n+1} &=& u_n + \dt \sum_{j=1}^m b_j \Bigl( -A v_{nj} + g(t_n+c_j \dt) \Bigr).
\end{eqnarray}
The method is determined by its coefficients $a_{ij}, b_j, c_i$
($i,j=1,\dots, m$).
We denote the Runge-Kutta matrix by
$\A=(a_{ij})$ and the row vector of the weights by $b^T=(b_j)$.
Eliminating the internal stages $v_{ni}$ results in
\begin{equation}\label{u-rec}
u_{n+1} = r(-hA) u_n + h \sum_{i=1}^m q_i(-hA)\, g(t_n+c_i \dt)\, ,
\quad\  n\ge 0,
\end{equation}
where the rational approximation $r(z)$ to $e^z$ is  defined by
\begin{equation}\label{rz}
r(z) = 1 + z b^T(I-z\A)^{-1}\1
\end{equation}
with $\1=(1,\dots,1)^T$, and where the rational functions~$q_i(z)$
are the entries of the row vector\footnote{
Instead of taking $r(z)$ and $q_i(z)$ as rational functions
originating from a Runge-Kutta method,
another suitable choice would be $r(z)=e^z$ and 
$q_i(z)= \int_0^1 e^{(1-\theta)z} \,\ell_i(\theta)\, d\theta$, where
$\ell_i$ is the $i$th Lagrange polynomial corresponding to the 
Gauss nodes~$c_j$.
This could be used similarly in the algorithm below.} 
\begin{equation}\label{sz}
q(z) = \bigl(q_1(z),\dots,q_m(z)\bigr) = b^T (I-z\A)^{-1}.
\end{equation}
We assume that the eigenvalues of the Runge-Kutta matrix $\A$ have
positive real part, and that the method is L-stable, i.e.,
\begin{equation}\label{L-stable}
|r(z)|\le 1 \quad\hbox{for \  Re\,}z\le 0, \qquad\hbox{and}\qquad
r(\infty)=0.
\end{equation}
These conditions are in particular satisfied by the Radau IIA
family of Runge-Kutta methods \cite{HaW96}.

The discrete analogue of the 
variation-of-constants formula
$$
u(t) = e^{-tA}u_0 + \int_0^t e^{-(t-\tau)A} \, g(\tau)\, d\tau
$$
is obtained by solving the recurrence relation (\ref{u-rec}). 
With the column vector 
$g_j=\bigl( g(t_j+c_i\dt) \bigr)_{i=1}^m$, this becomes
\begin{equation}\label{dvc}
u_{n} = r(-hA)^{n} u_0 +
h \sum_{j=0}^{n-1}  r(-hA)^{n-1-j}\,q(-hA)  \,g_j \, , \quad\  n\ge 1.
\end{equation}

\section{Discretization of the contour integrals}
\label{sec:contour}
\label{Sec.Contour}
We now discretize the Cauchy integral representation
\begin{equation}\label{contour-int}
r(-hA)^{n}q(-hA) = \frac1{2\pi {\rm i}}\int_\Gamma 
(\lambda+A)^{-1} \, r(\dt\lambda)^n\, q(h\lambda) \,  d\lambda
\end{equation}
along suitable contours $\Gamma$ in the resolvent set of $-A$.
The numerical integration in (\ref{contour-int})
is done by applying the
trapezoidal rule with equidistant steps to a parametrization of a
hyperbola~\cite{LoP04}.  
With one contour and one set of quadrature points on this contour,
we do not have a uniformly good approximation for all $n=0,\dots,N$,
but we can instead obtain a uniform approximation  
locally on a sequence of geometrically growing
intervals
\begin{equation}\label{intervals}
  I_\ell = [B^{\ell-1}h, B^{\ell}h),\qquad \ell\ge 1,
%  I_\ell = (B^{\ell-1}h, B^{\ell}h], % cf next section
\end{equation}
where the base $B>1$ is an integer, e.g.,  $B=10$.
For $nh \in I_\ell$ we approximate the contour integrals (\ref{contour-int}) as
\begin{eqnarray}\label{eq:num-int}
&& r(-hA)^{n}q(-hA)  \\
&& \qquad \qquad \approx \!
  \sum_{k=-K}^{K} w_k^{(\ell)}  (\lambda_k^{(\ell)}+A)^{-1}
r(\dt\lambda_k^{(\ell)})^n\, q(h\lambda_k^{(\ell)}) \, ,
  \ \ \ nh\in I_\ell, \nonumber
\end{eqnarray}
with the quadrature points $\lambda_k^{(\ell)}$ lying on 
a hyperbola $\Gamma_\ell$ and with the corresponding weights $w_k^{(\ell)}$.
The number of quadrature points on $\Gamma_\ell$, $2K+1$,
is chosen independent
of $\ell$. 
The contour $\Gamma_\ell$ is chosen as a hyperbola  given by
%(omitting $\ell$ in the notation)
\begin{eqnarray}
  \Rset \to \Gamma_\ell &:& \ \theta \mapsto \gamma_\ell(\theta) =
  \mu_\ell \,(1 - \sin(\alpha+i\theta)) + \sigma
  \label{hyp-param}
\end{eqnarray}
with an $\ell$-dependent parameter $\mu_\ell>0$. 
The angle $\alpha$ satisfies
$0<\alpha<\frac\pi 2 - \varphi$ with $\varphi$ of (\ref{sector}),
and $\sigma$ is the shift in (\ref{sector}).
The weights and quadrature points in (\ref{eq:num-int}) are given by
$$
w_k^{(\ell)} = \frac{i \tau }{2 \pi}\: \gamma_\ell'(\theta_k)~, \quad
\lambda_k^{(\ell)} = \gamma_\ell(\theta_k) \quad\ \mbox{ with }\quad
\theta_k = k\tau~,
$$
where $\tau$ is a step length parameter that can be chosen independent
of $\ell$.

The following bound of the necessary number of quadrature points 
is a consequence of the error analysis in
Section~\ref{sec.error}.

\pagebreak[3]

\begin{theorem}\label{thm:K}
In $(\ref{eq:num-int})$,
a quadrature error  bounded in norm by $\varepsilon$ for $nh\in I_\ell$
is obtained with 
$$
\textstyle K=O(\Log\frac1\varepsilon)\, .
$$ 
This holds
for $n \ge c \,\log (1/\varepsilon)$, with some constant $c>0$.
The required number $K$ is independent of $\ell$ and of $n$ and $h\le h_0$ 
with $nh\le T$.
For $\sigma\le 0$, $K$ is also independent of the length $T$ 
of the time interval. $K$~depends on the angle $\varphi$, the bound $M$
and the shift $\sigma$
in $(\ref{sector})$, but is otherwise independent of $A$.
\end{theorem}

The approximation is, however, poor for the first few $n$; 
cf.\:also~\cite{SchLL}.

Concerning the choice of parameters we  remark that the above asymptotic
bound for $K$ is obtained with 
$1/\tau$ proportional to $\log (1/\varepsilon)$ and 
with the parameter $\mu_\ell$
for the contour $\Gamma_\ell$  chosen such that 
$\mu_\ell B^\ell h = c_1\log (1/\varepsilon)$ with 
$c_1$ independent of $\ell$ and $h$, e.g., with $c_1=1/4$.
Since perturbations in the terms of (\ref{eq:num-int}) can be magnified
with $r(h\kappa_\ell)^n \approx e^{\kappa_\ell nh}$ with 
$\kappa_\ell= \mu_\ell(1-\sin\alpha)+\sigma$, the factor $c_1$
should not be chosen too large.
We refer to \cite{LoPSch} for an optimized strategy to choose the
parameters.

\pagebreak[3]

\section{The fast algorithm} 
\label{Sec.FA}
We start from the discrete variation-of-con\-stants formula (\ref{dvc})
for the Runge-Kutta approximation $u_N$ with a fixed $N$.
For the expression $r(hA)^N u_0$ we use the discretization of the
Cauchy integral like in the previous section and in fact similarly to
the approach of \cite{ShST00} for computing $\exp(-tA)u_0$.

The novel algorithm is concerned with the treatment of the inhomogeneity.
For a fixed step number~$N$  and a given base $B$ we split the
sum in~(\ref{dvc})
into $L$ sums, where $L$ is the smallest
integer such that  %$N \le B^{L}+1$: % 
$N \le B^{L}$:
$$
u_N = u_N^{(0)} + \dots + u_N^{(L)}
$$
with $u_N^{(0)} = hq(-hA) g_{N-1}$ and
$$
 u_N^{(\ell)} =
\dt\sum_{(N-1-j)h \in I_\ell} 
 r(-hA)^{N-1-j}\,q(-hA)  \,g_j
$$
for $\ell\ge 1$. %, and $u_N^{(0)} = hq(-hA) g_{N-1}$. 
On inserting the integral representation (\ref{contour-int})   
we obtain, with $n_\ell=N-B^\ell$ for $0\le \ell \le L-1$ 
%and with $n_0=N-1$ 
and~$n_L=0$,
$$
 u_N^{(\ell)}  = \dt\sum_{j=n_\ell}^{n_{\ell-1}-1}
 \frac 1{2\pi {\rm i}}\int_{\Gamma_\ell}
 (\lambda+A)^{-1} \, r(\dt\lambda)^{N-1-j}\, q(\dt\lambda) \,g_j \, d\lambda.
$$
The integral is discretized with the quadrature formula
of Section~\ref{sec:contour}: we approximate
$u_N^{(\ell)}$ by $U_N^{(\ell)}$ given as
\begin{eqnarray*}
U_N^{(\ell)} \ & = &
 \dt\sum_{j=n_\ell}^{n_{\ell-1}-1}
 \sum_{k=-K}^K w_k^{(\ell)} \, (\lambda_k^{(\ell)}+A)^{-1} \,
  r(\dt\lambda_k^{(\ell)})^{N-1-j} \, q(\dt\lambda_k^{(\ell)}) \,g_j
\\
&=&  \sum_{k=-K}^K w_k^{(\ell)} \,
  r(\dt\lambda_k^{(\ell)})^{N-n_{\ell-1}} \, (\lambda_k^{(\ell)}+A)^{-1} \,
y_k^{(\ell)},
\end{eqnarray*}
where
$$
y_k^{(\ell)} =
\dt\sum_{j=n_\ell}^{n_{\ell-1}-1} r(\dt\lambda_k^{(\ell)})^{n_{\ell-1}-1-j} 
\, q(\dt\lambda_k^{(\ell)}) \, g_j.
$$
Comparing this formula with (\ref{dvc}), we see that $y_k^{(\ell)}$ 
is the Runge-Kutta approximation to the solution at time
$t=n_{\ell-1}\dt$ of the linear initial-value problem
\begin{equation}\label{ivp-l}
y'(t) = \lambda_k^{(\ell)} y(t) + g(t), \qquad y(n_\ell \dt) =0,
\end{equation}
and hence $y_k^{(\ell)}$
is computed by Runge-Kutta time-stepping on (\ref{ivp-l}), using
(\ref{u-rec}) with the scalar $h\lambda_k^{(\ell)}$ in place of 
the operator $-hA$. With the solutions $x_k^{(\ell)}$ of the 
linear systems of equations
\begin{equation}\label{lin-system}
(\lambda_k^{(\ell)}+A) \, x_k^{(\ell)} = 
y_k^{(\ell)} ,
\end{equation}
 we  obtain $U_N^{(\ell)}$ as the linear combination
\begin{equation}\label{u-approx}
U_N^{(\ell)} =  \sum_{k=-K}^K c_k^{(\ell)} \, x_k^{(\ell)}
\qquad\hbox{with}\quad
c_k^{(\ell)}=\ w_k^{(\ell)} \,r(h\lambda_k^{(\ell)})^{B^{\ell-1}}.
\end{equation}  
There are only 
$(K+1)L$ linear systems (\ref{lin-system}) 
to be solved, for $k=0,\dots,K$ and $\ell\le L$.
(Since the quadrature points %$\lambda_k^{(\ell)}$
 lie symmetric with respect to the real axis,  
only the sum of the real parts of half the terms in (\ref{u-approx})
needs to be computed when approximating solutions with real components.)
We recall $L-1\le \log_B N$ and $K=O(\Log (1/\varepsilon))$, where
$\varepsilon$ is the accuracy requirement in the discretization
of the contour integrals. Note that the only approximation made 
in the computation of $U_N^{(\ell)}$,
is the discretization of the
contour integrals. 

Because of the poor approximation of the  contour integral (\ref{contour-int})
for small~$n$, we  evaluate
 $U_N^{(0)}+ U_N^{(1)}$ by $B$ direct Runge-Kutta steps up to
time $t=N\dt$ for the initial value problem
\begin{equation}\label{direct}
v'(t) + Av(t) = g(t),  \qquad  v((N-B)\dt)=0.
\end{equation}
This requires the solution of another $mB$ linear systems with matrices
of the form $(\lambda+A)$. For small values of $B$ or stringent
accuracy requirements, we take
$B^2$ direct Runge-Kutta steps to compute ${u_N^{(0)}+ u_N^{(1)}+u_N^{(2)}}$.
(Asymptotically, we need to take $O(\log (1/\varepsilon))$ direct steps
according to Theorem~1.)

Finally we sum up the $U_N^{(\ell)}$ to obtain
\begin{equation}\label{u-sum}
U_N = U_N^{(0)} + \dots + U_N^{(L)}
\end{equation}
as the approximation to $u_N$.
The fast algorithm thus consists of doing the steps 
(\ref{ivp-l})--(\ref{u-sum})
in the given order.

\begin{remark}
The algorithm extends to differential equations with a positive definite
mass matrix $M$,
\begin{equation}\label{ivp-mass}
M u'(t) + A u(t) = g(t), \qquad u(0)=u_0,
\end{equation}
which is transformed to a system $\~u'(t) + \~A \~u(t) = \~g(t)$ for 
$\~u(t)= M^{1/2} u(t)$
with $\~A= M^{-1/2} A M^{-1/2}$ and $\~g(t) = M^{-1/2}g(t)$.
Applying {\it formally} the above algorithm to the transformed system
and then transforming back yields again (\ref{u-approx}), where now
$x_k^{(\ell)}$ is the solution of the linear system
\begin{equation}\label{lin-system-mass}
(\lambda_k^{(\ell)}M+A) \, x_k^{(\ell)} = y_k^{(\ell)} ,
\end{equation}
and $y_k^{(\ell)}$ is the Runge-Kutta approximation at $t=n_{\ell-1}\dt$
of the initial value problem (\ref{ivp-l}) with the untransformed
inhomogeneity $g(t)$.
\end{remark}

\begin{remark} We have formulated the algorithm for a constant time step
size~$h$, but this is not essential. The algorithm is readily extended
to accommodate variable step sizes, with
the same step size sequence for all $k$ in (\ref{ivp-l}), chosen adaptively 
according to the behaviour of the inhomogeneity $g(t)$.
Adaptivity in space can be  used in solving the linear systems
(\ref{lin-system}), choosing the spatial mesh according to
the behaviour of the right-hand sides $y_k^{(\ell)}$ and the operator $A$. 
Note that in a hierarchical basis representation, adding a mesh point just
corresponds to  adding a scalar differential equation in (\ref{ivp-l}).
The details of
such an adaptive algorithm are beyond the scope of this paper.
\end{remark}

\section{Numerical experiment}
\label{sec:NumEx}

\begin{figure}[htbp]
  \centering
  \includegraphics[width=0.4\textwidth]{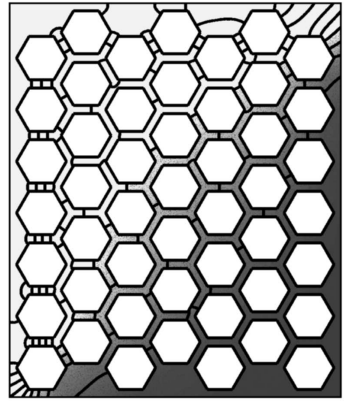}
  \caption{Domain for the heat equation, with isolines of the temperature distribution 
  at $t = 20$.}
  \label{fig:geo}
\end{figure}

We consider an initial-boundary value problem of 
the heat equation in two space dimensions for $u=u(x,t)$,
$$
\left\{\, \begin{array}{llll}
  \partial_{t} u(x,t)  &=& \Delta u (x,t),& \quad
   x\in \Omega, \, 0\le t \le T,
  \\[2pt]
  u(x,0) &=& 0 , & \quad  x\in \Omega,
  \\[2pt]
  \partial_{\nu} u (x,t) &=& 0, &  \quad  x\in \Gamma_{int}, \, 0\le t \le T,
  \\[2pt]
  \partial_{\nu} u (x,t)&=& \beta(x,t) - \rho (u(x,t)-u_{out}) ,& \quad
     x\in \Gamma_{out}, \, 0\le t \le T,
\end{array}
\right.
$$
on a wire-fence like structure (rectangle of size $10.65 \times 12.64$ 
with hexagonal holes, each hole with radius 0.8),
see~Figure~\ref{fig:geo}. Here $\Gamma_{int}$ denotes 
the boundary of the holes, and
$\Gamma_{out}$ is the boundary of the rectangle. 
In the example we set the heat flux $\beta = 5 \sin^{2}(t)$ 
on the upper and left boundary of the rectangle and
$\beta=0$ on the lower and right boundary,
and the convective heat flux to $\rho (u-u_{out})$,
with the ambient temperature $u_{out} = 0$ and the coefficient of
surface heat transfer $\rho = 0.5$, cf.~the introduction in~\cite{LMTS96}.
Space is discretized using linear finite elements on a triangular mesh, with
$27346$ vertices and $50368$ triangles. Triangulation is done using the tool 
Triangle~\cite{Shewchuck96b}. 

The finite element equations are of the form (\ref{ivp-mass}),
where $M$ is the standard mass matrix containing the $L^2$ inner products of the
nodal basis functions $\varphi_i$. The stiffness matrix is the
sum $A = A_{0} + \rho M_{b}$ with
$$
A_{0}\big|_{ij} = \int_{\Omega} \nabla \varphi_{i} \nabla \varphi_{j}\, dx~,
\quad 
M_{b}\big|_{ij} = \int_{\Gamma_{out}} \varphi_{i} \varphi_{j}\, d\sigma ~.
$$  
The inhomogeneity $g(t)$ is given by
$$
g_i(t) = 
\int_{\Gamma_{out}} (\beta(x,t) + \rho u_{out}) \varphi_{i} \,d\sigma(x) .
$$
The algorithm takes into account that $g(t)$ has nonzero entries only
along the
outer boundary $\Gamma_{out}$, so that effectively $g(t)$ is a vector whose
dimension is the number of degrees of freedom 
on the outer boundary -- in this example $776$. 
The differential equations (\ref{ivp-l})
need to be integrated only for this reduced dimension, since they have no
coupling between the components.

We have used the $2$- and $3$-stage Radau IIA methods 
(of orders~3 and 5, respectively) for
time discretization in 
our numerical experiments.

In the fast algorithm we set $B = 5$ and $K=15$ and, from the experience
of \cite{LoPSch,SchLL}, we choose the angle in the hyperbola as
$\alpha= \pi/4 $, the parameter
$\mu_{\ell} = 3/ (\dt B^\ell)$
and the parameter $\tau = 5/K$.
This choice of parameters leads to a deviation of the order~$10^{-6}$ 
from the Runge-Kutta approximation
at time $t=20$.

\begin{figure}[tbp]
  \centering
 \includegraphics[width=0.49\textwidth]{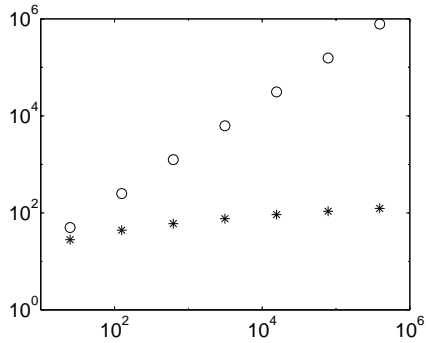}
  \caption{Number  of solves of linear systems versus~step number: 
  direct time-stepping~$(\circ)$
  and fast algorithm~$(*)$.}
  \label{fig:lscount}
\end{figure}

The two-dimensional example above is still small enough that a direct solution
of the linear systems using sparse solvers is reasonable.
A direct implementation of the $m$-stage Radau~IIA method (cf.~\cite{HaW96})
requires only $m$ sparse LU factorizations, computed at the beginning
of the integration, followed by $mN$ substitutions.
On the other hand, for the algorithm presented here we need to solve
$(K+1)(L-1)$ linear systems with matrices $\lambda M +A$
for as many different values of~$\lambda$, and the $mB$ linear systems for 
the $B$
direct steps. 
Especially with a diagonal, lumped mass matrix $M=DD^T$, this work can be reduced
by a similarity transform taking $D^{-1}AD^{-T}$ to tridiagonal (or
Hessenberg) form $T$, but
exploiting sparsity here becomes an issue; see~\cite{Cav94,Nik00}. 
The resulting
linear systems with $\lambda I+T$ are then inexpensive to solve.
Even without using such a transform, the fast algorithm eventually
overtakes the standard algorithm for sufficiently large step
numbers $N$, in the present example for $N\approx 1000$.
Much earlier and larger relative gains arise 
when iterative solvers are used for the linear systems 
in both algorithms, as is clear from the linear systems count
in Figure~\ref{fig:lscount}.

\section{Error analysis}
\label{sec.error}
Our analysis relies on the good behaviour of the trapezoidal rule for
certain holomorphic integrands \cite{LoP04,Ste,Sten}. 
Following the ideas in \cite{LoP04},  
we consider the continuation of the para\-metrization (\ref{hyp-param}) 
to the conformal mapping
\begin{equation}\label{gamma-w}
\gamma(w)=\mu \, (1-\sin(\alpha+iw)).
\end{equation} 
(For ease of presentation we set $\sigma = 0$ in (\ref{sector}).)
This conformal mapping
transforms each horizontal straight line
$$
\imag w = y,\qquad -d \le y \le d,
$$
with $0<\alpha -d < \alpha + d < \frac\pi 2$,
into the left branch of the hyperbola
$$
\lambda\in\bC:\
\bigg(\frac{\real \lambda -\mu}{\mu\sin(\alpha-y)} \bigg)^2
- \bigg(\frac{\imag \lambda}{\mu\cos (\alpha-y)}\bigg)^2=1,$$ i.e.,
the left branch of the hyperbola with center at $(\mu,0)$,
foci at $(0,0),\ (2\mu,0)$ and with asymptotes forming angles
$\pm [\pi/2-(\alpha-y)]$ with the real axis. Therefore, $\gamma$
transforms the horizontal strip
$$
D_d=\{w\in \bC : |\imag w|\leq d\}
$$
into the region $\Omega=\gamma(D_d)$
limited by the left branches corresponding to
$y=\pm d$. To indicate the dependence on the parameter $\mu$ of
(\ref{gamma-w}), we write $\Omega=\Omega_\mu$. We note that 
$\lambda\in\Omega_\mu$ if and only if
$h\lambda\in\Omega_{h\mu}$ for any $h>0$, so that
$$
h\Omega_\mu = \Omega_{h\mu}.
$$
Because of (\ref{sector}), henceforth we will 
assume that $\alpha >0 $ and $d >0$ satisfy 
$0<\alpha -d < \alpha + d < \frac\pi 2 -\varphi$. 
Under these conditions, all the hyperbolas we are considering lie outside 
the spectrum
of $-A$. 

After parametrizing (\ref{contour-int})
via $\gamma$, we get
$$
r(-hA)^n \, q(-hA) = \int_{-\infty}^{+\infty}G_{h,n}(x)\, dx,
$$
where $G_{h,n}(w)$ is given, for $w\in D_d$, by 
\begin{equation}\label{hyperbolaintegrand}
G_{h,n}(w)=\frac{1}{2\pi i} \Bigl(\gamma(w)+A\Bigr)^{-1}\,
r(h \gamma(w))^n \, q(h \gamma(w))\,\gamma'(w).
\end{equation}
For an integrable mapping $G : \Rset \to X$, $K \ge 1$ and $\tau >0$, set
\begin{equation}
\label{quadratureerror}
E_{\tau,K}(G) =\int_{-\infty}^{+\infty}G(x)\, dx-\tau
\sum_{k=-K}^{K}G(k\tau),
\end{equation}
i.e., $E_{\tau,K}(G)$ stands for the quadrature error of the truncated
trapezoidal rule for the integral of $G$.
Our goal is precisely to estimate $E_{\tau,K}(G_{h,n})$.
To this end we first consider the behaviour of $G_{h,n}$ on $D_d$. We need the 
following lemma whose elementary proof is omitted.

\begin{lemma}\label{estimaterational}
Let $r(z)$ be a rational function with $r(0)=1,\
r'(0)=1$ which satisfies the L-stability condition (\ref{L-stable}). 
Then, there exist $\rho>0$ and
%$0< b <\sin(\alpha-d)$ 
$b>0$ such that
\begin{equation}\label{Radaukernelestimate}
|r(z)|\leq  \frac{e^{2\delta}}{1+b\,|z|} ,
\qquad \hbox{for }\ z \in \Omega_{\delta}\ \hbox{with }\ 0 < \delta \le \rho.
\end{equation}
\end{lemma}

Now, from the sectorial condition (\ref{sector}) on $A$
and Lemma~\ref{estimaterational} with 
${\delta=h\mu\le\rho}$, 
we obtain
\begin{equation}\label{cotaGhn}
\|G_{h,n}(x+iy)\| \leq C_0 \,\frac{e^{2\mu hn}}{(1+bh\mu
(\cosh x - \sin(\alpha-y)))^{n}}
\end{equation}
for $ x\in\bR $ and $|y|\leq d$ (recall that $0<\alpha-d < 
\alpha+d < \frac\pi2 -\varphi$),
where $C_0$ is the constant given by
$$
C_0=\frac{M}{2\pi}
\sqrt{\frac{1+\sin(\alpha+d)}{1-\sin(\alpha+d)}}\:\; 
\max_{z\in\Omega_\rho} \| q(z) \|.
$$
Finally, the above bound (\ref{cotaGhn}), the elementary
inequality
$$
1+c-s \geq (1-s)(1+c), \qquad c,s > 0,
$$
and the bound 1 for the sine yield, for $|y|\le d$ and $t=nh$, 
\begin{eqnarray}
\label{finalGhn}
&&\|G_{h,n}(x+iy)\|  %\\[6pt] && \hspace{48pt} 
\le \,
\frac{C_0\,e^{2\mu t}}{(1-b\mu t/n)^n}\,\Big( 1+ \frac{b\mu t}{n} \cosh x
\Big)^{-n}. 
%\nonumber
\end{eqnarray}

Next, to estimate $E_{\tau,K}(G_{h,n})$, we are going to use an approach
similar to the one in  \cite{LoP04,Ste,Sten}. 
We denote by $S(D_d,X)$ the class  formed by all the continuous 
mappings $G: D_d\to X$ (for a complex Banach space $X$, here a space of matrices) 
holomorphic on the interior of the strip $D_d$, and satisfying
the following two conditions:

\begin{equation}\label{pp1} \int_{-d}^d
\|G(x+iy)\|\,dy \to 0, \quad {\rm as}\ |x|\to +\infty,
\end{equation}
\begin{equation}\label{pp2}
N( G ,D_d)  :=  
\int_{-\infty}^{+\infty}\{\|G(x+id)\| +
\|G(x-id)\|\}\,dx <+\infty.
\end{equation}

Given $G \in S(D_d,X)$, it turns out, assuming that $G$ has a fast decay 
at $\infty$, that $E_{\tau,K}(G)$
becomes very small as $K \to +\infty$ if
$\tau$ is properly tuned (see \cite{LoP04,Ste,Sten} for various situations). 
In Theorem~\ref{errorestimate} we assume that $G$ exhibits the kind
of decay of $G_{h,n}$ in (\ref{finalGhn}) and this theorem will directly
provide the estimate for $E_{\tau,K}(G_{h,n})$ we are looking for.

\begin{theorem}\label{errorestimate}
Assume that $G\in S(D_d,X)$ for some $d>0$, and that there exist
$C, a  >0$ and $n \ge 1$ such that
\begin{equation}\label{desigualdadpowercosh}
\|G(x)\| \leq C \Big( 1+ \frac{a }{n}\cosh x\Big)^{-n},
\qquad x \in \bR.
\end{equation}
Then, for $\tau >0,\ K\geq 1$, there holds
\begin{eqnarray*}
\|E_{\tau,K}(G)\| &\le& \frac{N( G ,D_d)}{e^{2\pi d/\tau}-1} 
\\[1mm]
&&
+ C \, \left(  \Lfunc( a )\, e^{- a  \cosh(K\tau)/2}+ 
\five \Big( 1 + \frac{a }n \cosh K\tau \Big)^{-(n-1)} \right),
\end{eqnarray*}
with $\Lfunc( a )=2+ |\log (1-e^{-a /2})|$.
\end{theorem}

Notice that $\Lfunc$ is decreasing, $\Lfunc( a)\to 2$ as $a\to +\infty$ and
$\Lfunc( a)\sim |\log a|$ as $a\to 0^{+}$.

\begin{proof} Set
$$
E_{\tau,\infty}(G)=\int_{-\infty}^{+\infty}G(x)\,dx - \tau 
\sum_{k=-\infty}^{\infty}G(k\tau),\qquad \tau>0. 
$$
For fixed $K\geq 1$, it is clear that
$$
\|E_{\tau,K}(G)\| \leq \|E_{\tau,\infty}(G)\| + \tau \sum_{|k|\geq
K+1} \|G(k \tau)\|.
$$
On the one hand, by Theorem 4.1 in \cite{Ste} (see also \cite{Sten}), we
have
$$
\|E_{\tau,\infty}(G)\| \leq \frac{N(G,D_d)}{e^{2\pi d/\tau}-1}.
$$
On the other hand,
\begin{eqnarray*}
\tau \sum_{|k|\geq K+1} \|G(k \tau)\| 
%&=& \tau\sum_{k=K+1}^{+\infty} (\|G(k \tau)\| + \|G(-k \tau)\|) 
%\\[5pt]
&\leq& 2C \tau \sum_{k=K+1}^{+\infty} \Big( 1+
\frac{a }{n}\cosh k\tau \Big)^{-n} \\ [5pt] &\leq& 2C
\int_{K\tau}^{+\infty}  \Big( 1+ \frac{a }{n}\cosh x
\Big)^{-n} \,dx.
\end{eqnarray*}
The proof of the theorem is now 
completed by applying the following lemma. $\quad\ \Box$
\end{proof}

\begin{lemma}\label{integralpowercoshestimate}
For $R\geq 0$, $a  > 0$ and $n \ge 1$ there holds
$$
\int_{R}^{+\infty} \Big(1+ \frac{a }{n}\cosh x
\Big)^{-n}\, dx \leq \Lfunc( a ) \, e^{-a \cosh R/2}+
\five  \Big(1+\frac{a }{n}\cosh R \Big)^{-(n-1)}.
$$
\end{lemma}

\begin{proof}
The change of variables $u=\cosh x$ shows that
$$
\int_{R}^{+\infty} \Big(1+ \frac{a }{n}\cosh x
\Big)^{-n}\, dx = \int_{\cosh R}^{+\infty} \Big(1+ \frac{a }{n} u
\Big)^{-n} \frac{du}{\sqrt{u^2-1}} \,.
$$
Set $\beta=\max\{\cosh R,n/a \}$. Then, from the estimates in
\cite{LoP04} and the elementary inequality
\begin{equation}
\label{elineq}
(1+ y/n)^{-n}  \leq   e^{-y/2},  \qquad \mbox{for }\ 0\le y \le n
\end{equation}
it turns out that
\begin{eqnarray*}
\int_{\cosh R}^{\beta} \Big(1+ \frac{a }{n} u \Big)^{-n} \frac{du}{\sqrt{u^2-1}} & \le &
\int_{\cosh R}^{\beta} e^{-a u/2}  \frac{du}{\sqrt{u^2-1}}  \\[4pt]
& \le  & \int_{R}^{+\infty} e^{-a  \cosh x /2}\, dx \\[4pt]
&\le & \Lfunc( a )\, e^{-a \cosh R/2} \, .
\end{eqnarray*}
Moreover,
\begin{eqnarray*}
&&\int_\beta^{+\infty} \Big(1+ \frac{a }{n} u \Big)^{-n}
\frac{du}{\sqrt{u^2-1}}  
\\
&&\qquad\qquad \leq
\Big(1+ \frac{a }{n} \cosh R \Big)^{-(n-1)} 
\int_\beta^{+\infty} \Big(1+ \frac{a}{n} u
\Big)^{-1} \frac{du}{\sqrt{u^2-1}}\,.
\end{eqnarray*}
Now, since $\beta\ge \max\{1,n/a \}$, the result follows from
the observation that for both $n/a \ge 1$ and $n/a \le 1$ we have
$$
\int_{\max\{1,n/a \} }^{+\infty}  \Big(1+ \frac{a}{n} u \Big)^{-1}
\frac{du}{\sqrt{u^2-1}}\leq 
\int_1^{+\infty} (1+v)^{-1} \, \frac{dv}{\sqrt{v^2-1}}= 1.
\quad\ \Box
$$
\end{proof}

We apply  Theorem~\ref{errorestimate} to $G_{h,n}$. First of all, notice that by (\ref{finalGhn}) 
it is clear that $G_{h,n}$ satisties (\ref{pp1}). 
Moreover, by Lemma~\ref{integralpowercoshestimate}, we have 
\begin{eqnarray}
\label{NGhn}
N(G_{h,n},D_d) & \le & \frac{4C_0\, e^{2\mu t}}{(1-b\mu t/n)^{n}} \\[5pt]
&& \hspace{24pt} \times
\bigg( \Lfunc(b\mu t)\,e^{-b\mu t/2} + \Big(1+ \frac{b\mu t}{n}\Big)^{-(n-1)} \bigg), \nonumber
\end{eqnarray}
and conclude that $G_{h,n} \in S(D_d,X)$.
Then, in view of (\ref{finalGhn}) and (\ref{NGhn}), Theorem~\ref{errorestimate} 
yields directly
\begin{eqnarray*}
\|E_{\tau,K}(G_{h,n})\| &\leq& \frac{4C_0 \, e^{2\mu t}}
{(1-b \mu t/n)^{n}}
%\\[5pt]\nonumber &\times& 
\Bigg(  \frac{\Lfunc( b\mu t)\,e^{-b\mu t/2}
+\five (1+b\mu t/n)^{-(n-1)}}{e^{2\pi d/\tau}-1}
\\[5pt]\nonumber
&+& \Lfunc( b\mu t) \, e^{-b\mu t\cosh(K\tau)/2} \!
%\\[5pt]\nonumber &+& 
+\!\Bigl(1+\frac{b\mu t}n \cosh(K\tau)\Bigr)^{-(n-1)}
\Bigg).
\end{eqnarray*}
A simplified version of this estimate 
is obtained by using the elementary inequalities (\ref{elineq}) and
$$
\begin{array}{cccl}
(1-y/n)^{-n} & \le & e^{2y}, & \qquad \mbox{for } 0\le y\le n/2, \\[1mm]
\Lfunc(y)& \le & 3, & \qquad  \mbox{for } y \ge1.
\end{array}
$$
Setting
$$
\textstyle
C=20\, C_0,\ \, a_0=2+\frac32 b,\ a_1=2+2b, \ a_2 = \frac12 b,
$$
with $b$ of Lemma~\ref{estimaterational} as before, we can
summarize the final result in the following theorem.

\begin{theorem}\label{thm:err}
The quadrature error $(\ref{quadratureerror})$ for $G_{h,n}$ of 
$\,(\ref{hyperbolaintegrand})$ with $(\ref{gamma-w})$ satisfies, 
for $t=nh$ and if $n/2 \ge b\mu t\ge 1$,
\begin{eqnarray*}
\|E_{\tau,K}(G_{h,n})\| &\leq& C
\Bigg(  \frac{e^{a_0\mu t}}{e^{2\pi d/\tau}-1}
+ e^{(a_1-a_2\cosh(K\tau))\mu t}
\\
&&
\quad
+\ e^{a_1\mu t}\Bigl(1+\frac{b\mu t}n \cosh(K\tau)\Bigr)^{-(n-1)}\Bigg).
\end{eqnarray*}
\end{theorem}

The first term in the error bound becomes 
$O(\varepsilon)$ if $\tau$ is chosen so small that
$a_0\mu t - 2\pi d/\tau \le \log \varepsilon$, which requires
an asymptotic proportionality
$$
\frac 1\tau \sim \log\frac1\varepsilon + \mu t.
$$
For $\mu$ chosen such that
$$
\frac{c_1}B \log\frac1\varepsilon \le \mu t \le c_1 \log\frac1\varepsilon
$$ 
with an arbitrary positive constant $c_1$ and with $B> 1$, we obtain that the
second term  is  $O(\varepsilon)$ if
$a_1 -  a_2  \cosh(K\tau) \le -B/c_1$, i.e., with
$$
\cosh(K\tau) = c_2
$$
for a sufficiently large constant $c_2$. With the above choice of $\tau$,
this yields
$$ 
K\sim \log\frac1\varepsilon.
$$
The third term then becomes smaller than $\varepsilon$ for
$$
n \ge c \, \log \frac1\varepsilon
$$
with a sufficiently large constant $c$.
Taken together, these estimates prove Theorem~\ref{thm:K}.

\vskip 12pt
\noindent 
{\it Acknowledgements.} The research of the first and third authors has 
been supported
by DGI-MCYT under project MTM2004-07194 cofinanced by FEDER funds. 
The research of the 
second author has been
supported by DFG, SFB 382. The reserach of the fourth author has been supported
by the DFG Research Center {\sc Matheon} 
\lq \lq Mathematics for key technologies" in 
Berlin.

\pagebreak[3]

\bibliographystyle{amsplain}

\begin{thebibliography}{10}

\bibitem{AsS94} A. Ashyralyev and P. Sobolevskii, 
{Well-Posedness of Parabolic Difference Equations}. Birkh\"auser,
Basel, 1994.

\bibitem{BaTW02} N. Y. Bakaev, V. Thom\'ee, and L. Wahlbin,
{Maximum-norm estimates for
resolvents of elliptic finite element operators}. Math. Comp. 72 (2002), 
1597--1610.

\bibitem{BrCT82} P. Brenner, M. Crouzeix, V. Thom\'ee,
{Single step methods for inhomogeneous linear differential equations
in Banach space}. RAIRO Mod\'el. Math. Anal. Num\'er. 16 (1982), 5--26.

\bibitem{Cav94} I.A. Cavers,
A hybrid tridiagonalization algorithm for symmetric sparse matrices.
SIAM J. Matrix Anal. Appl. 15 (1994), 1363--1380.

\bibitem{GaM} I.P. Gavrilyuk, V. Makarov,
Exponentially convergent algorithms for the operator exponential
with applications to inhomogeneous problems in Banach spaces.
Preprint, 2004.

\bibitem{HaW96} E. Hairer, G. Wanner, Solving Ordinary Differential
Equations. II. Stiff and Differential-Algebraic Problems.
Second edition. Springer, Berlin, 1996.

\bibitem{LMTS96} R. W. Lewis, K. Morgan, H.R. Thomas, K.N. Seetharamu,
The Finite Element Method in Heat Transfer Analysis. 
John Wiley \& Sons Ltd, Chichester, 1996.


\bibitem{LoP04} M. L\'opez-Fern\'andez, C. Palencia, 
On the numerical inversion of the Laplace
  transform of certain holomorphic mappings. 
  Appl. Numer. Math. 51 (2004), 289-303.

\bibitem{LoPSch} M. L\'opez-Fern\'andez, C. Palencia, A. Sch\"adle,
On the numerical inversion of the Laplace
  transform of certain holomorphic mappings, Addendum. (In preparation).

\bibitem{LuO93} C. Lubich, A. Ostermann,
Runge-Kutta methods for parabolic equations and convolution quadrature.
Math. Comput. 60 (1993), 105--131.

\bibitem{LuS02} C. Lubich, A. Sch\"adle,
Fast convolution for nonreflecting boundary conditions.
SIAM J. Sci. Comp. 24 (2002), 161--182.

\bibitem{Nik00} J.L. Nikolajsen,
An improved Laguerre eigensolver for unsymmetric matrices.
SIAM J. Sci. Comp. 22 (2000), 822--834.

%\bibitem{Riz95}
%M.~Rizzardi,
%A modification of {T}albot's method for the simultaneous
%approximation of several values of the inverse {L}aplace transform.
%ACM Transactions an Mathematical Software 21 (1995), 347--371.

\bibitem{SchLL}
A.~Sch\"adle, M.~L\'opez-Fern\'andez, C.~Lubich,
Fast and oblivious convolution quadrature. Preprint, 2005.

\bibitem{ShST00}
D. Sheen, I. H. Sloan, V. Thom\'ee, 
A parallel method for time-discretization of parabolic problems 
based on contour integral representation and quadrature.  
Math. Comp.  69  (2000),   177--195. 

\bibitem{ShST03}
D. Sheen, I. H. Sloan, V. Thom\'ee, 
A parallel method for time discretization of parabolic equations 
based on Laplace transformation and quadrature. 
IMA J. Numer. Anal. 23 (2003), 269--299. 

\bibitem{Shewchuck96b}
   J. R.~Shewchuk,
   Triangle:  {E}ngineering a {2D} {Q}uality {M}esh {G}enerator and
   {D}elaunay {T}riangulator, in
   Applied Computational Geometry:  Towards Geometric Engineering,
   Eds. M. C. Lin and D. Manocha,
   Lecture Notes in Computer Science
   1148, Springer, 1996, 203--222.

\bibitem{Ste}{F. Stenger}, {Approximations via Whittaker's
cardinal function}. J. Approx. Theory 17 (1976), 222--240.

\bibitem{Sten}{F. Stenger}, {Numerical methods based on Whittaker cardinal,
 or sinc functions}. SIAM Review 23 (1981), 165--224.

\bibitem{Tal79}
A.~Talbot, The accurate numerical inversion of Laplace transforms.
 J. Inst. Math. Appl. 23 (1979), 97--120.

\end{thebibliography}

\end{document}